# Time-bandwidth Study of Non-classically Damped, Linear, Time-invariant Coupled Oscillators with Closely Spaced Modes

Luis M. Baldelomar Pinto[1], Alireza Mojahed[2], Sobhan Mohammadi[3], Keegan J. Moore[3], Lawrence A. Bergman[1], Alexander F. Vakakis[1]

[1] University of Illinois Urbana-Champaign, Urbana, IL 61801
[2] 150 Glover Avenue, Norwalk, CT 06850
[3] Georgia Institute of Technology, Atlanta, GA 30332

**Abstract**

In dynamics and vibrations, the concept of bandwidth for linear time-invariant systems is widely recognized as a measure of the dispersion of frequency content around resonance. Similarly, the time constant is associated with the rate of energy decay in the time domain. Notably, the time-bandwidth limit for such systems is unity, indicating that achieving sharp frequency localization while simultaneously maintaining a slow energy decay is not feasible, nor is it possible to achieve a broad frequency spread while preserving a rapid energy decay. However, the time-bandwidth concept does not have a well-defined application to multi-degree of freedom (multi-DOF) systems characterized by strong modal interactions. This research aims to develop a comprehensive time and bandwidth concept for a linear two-DOF system with significant modal interactions. We focus on a non-classically damped system, which facilitates complex mode interactions, and we investigate how the definition of bandwidth and time constant can be applied to account for the slow dynamics observed in energy decay. By examining this system under various parameters, we gain insights into the energy decay behavior at specific time-bandwidth product regimes. Our analytical results are validated through experiments. Our findings elucidate the implications of the time-bandwidth product for a linear multi-DOF system's response and provide valuable insights into the influence of modal interactions on energy decay.

# 1. Introduction

For a single-degree-of-freedom and linear time-invariant system (single-DOF LTI), *bandwidth (BW)* is an inherent property which relates to the sharpness of its frequency spectral content. This indicates the range of frequencies at which its amplitude is above a determined threshold, typically 3 dB units below its maximum amplitude (Ewins, 2000; Yaghjian, 2005). BW is also important in the field of control systems, where it is often used as a performance metric; for instance, the BW of a passband filter indicates the range of frequencies that are allowed to pass through (Ewins, 2000; Nise, 2020). Similarly, the notion of time constant (or *energy storage time – EST*) is used for single-DOF LTI resonators as a metric of energy decay rate in the time domain (Capek et al., 2015). The combined product of these two parameters (BW and EST) is known as the *time-bandwidth product – TBP* (Tsakmakidis et al., 2017, 2020), and for single-DOF LTI resonators the TBP is equal to unity, which is referred to as the classical *time-bandwidth limit (TBL)*, an inherent fundamental property of this class of resonators, dictating the dual relation of energy dissipation (locality) in the frequency and time domains.

Recent studies have extended the concept of BW for dynamical systems. Mojahed et al. (2022a) presented a generalization of the rms BW definition based on earlier works for signal processing (Gabor, 1946; Starosielec and Hägele, 2014), to single-DOF nonlinear dynamical systems; to this end, they enveloped their velocities as a way to define measures of power dissipated by damping in the frequency domain and rate of dissipation in the time domain. Mojahed et al. (2022b) also performed an experiment proving how the classical TBL can be broken for single-DOF nonlinear resonators. Additionally, in (Chang et al., 2025) a BW definition for general LTI multi-DOF systems with classical (proportional) viscous damping distribution was developed, based on the superposition of the contributions of their modal energy decays. These early studies suggest rich unexplored territory towards understanding the concepts of BW and EST in general classes of multi-DOF LTI systems, especially if complex vibration modes due to non-proportional damping distribution are present (Caughey, 1960; Ewins, 2000). In a recent study (Zhang et al., 2026) the BW, EST and TBP were generalized to define nonlinear modal dissipation in multi-DOF nonlinear systems by generalizing the notion of spectral submanifolds (Haller and Ponsioen, 2016).

Accordingly, the aim of this study is to explore BW and EST in a linear two-DOF system with *closely spaced modes* and *non-classical damping*. Firstly, we discuss the definitions of BW and EST considered in this study in terms of an "effective oscillator" which is derived based on the total energy decay of the considered two-DOF system subject to impulsive excitation. In turn, the velocity envelope of the effective oscillator is used to describe the overall energy decay in the two-DOF system. Next, we present a numerical study of the BW, EST and time-bandwidth product (TBP) of this system using the derived envelope expression under different initial conditions and varying system parameters. We probe the energy response at specific system parameters to understand how our results relate to the dissipative properties of the two-DOF non-classically damped system, and show that *the combined effects of closely spaced modes and viscous damping non-proportionality result in breaking of the classical TBL*; the implications of this result in terms of the dissipative capacity and rate of energy dissipation of the two-DOF system are discussed. We also compare the resulting decaying energy responses to those of a single-DOF oscillator with equivalent BW or EST to gain further physical insights. Lastly, we present an experimental validation of the results, and provide concluding remarks synopsizing the main findings of this work.

## 2. Analysis
### 2.1. Two-DOF system with closely spaced modes and non-classical viscous damping

The system under consideration is a two-DOF discrete linear system with identical masses and equal grounding stiffnesses but different viscous dampers (Figure 1). Without loss of generality, throughout this work we normalize the masses $m_i$ and grounding stiffnesses $\beta_i$, according to $\beta_1 = \beta_2 = m_1 = m_2 = 1$, and assume *weak coupling and weak damping,* $\beta, \lambda_1, \lambda_2 \ll 1$, with the three parameters being of the same small order; moreover, the damping in the coupling element is eliminated, $\lambda = 0$. We note that the configuration symmetry of the system is perturbed only when $\lambda_2 \neq \lambda_1$, which renders this system *non-classically (non-proportionally) damped*. Unless otherwise stated, throughout this work we will assume (without loss of generality) that $\lambda_2 = 0.05 > \lambda_1 = 0.00125$, and denote the more lightly damped oscillator on the left as "oscillator 1", and the more heavily damped one on the right as "oscillator 2". Then, it turns out (Baldelomar Pinto et al., 2025) that the impulsive dynamics of this system are governed by a single parameter, namely the *damping non-proportionality to coupling stiffness ratio* $\gamma = \frac{4\beta}{\lambda_2 - \lambda_1}$. Indeed, for $\gamma < 2$, the impulsive responses possess two distinct dissipation rates and a single fast frequency; while for $\gamma > 2$, these responses have a single dissipation rate but two distinct fast frequencies. We will make use of this result for the rest of this work.

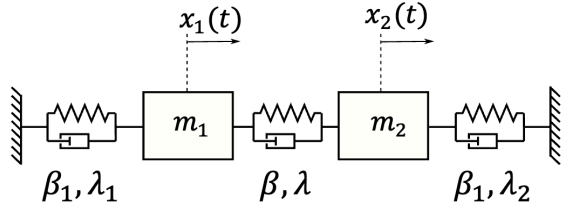

*Figure 1. Two-DOF system of coupled oscillators with non-classical damping, and normalized parameters $\beta_1 = \beta_2 = m_1 = m_2 = 1, \lambda = 0, \lambda_2 > \lambda_1$, and varying coupling stiffness $\beta$.*

An important feature of this system is the complexity of its modes due to the proximity of its vibration modes and non-classical damping distribution. This introduces modal interactions that enable energy exchanges not only between the two oscillators but also between the two vibration modes which become coupled. To visualize this, we characterize the complexity of the vibration modes through their modal phase collinearity (MPC), which is often used in modal analysis to determine whether a mode is real or complex-valued (Gres et al., 2021). We can calculate the MPC for this system by solving the eigenvalue problem of this system numerically, extracting their mode shapes, and inputting them to the following expression,

$$\text{MPC}(\underline{\varphi}) = \frac{(S_{xx} - S_{yy})^2 + 4S_{xy}^2}{(S_{xx} + S_{yy})^2} \quad (1)$$

where $S_{xx} = \Re\{\underline{\varphi}\}^T \Re\{\underline{\varphi}\}, S_{yy} = \Im\{\underline{\varphi}\}^T \Im\{\underline{\varphi}\}, S_{xy} = \Re\{\underline{\varphi}\}^T \Im\{\underline{\varphi}\}, \Re\{\bullet\}$ and $\Im\{\bullet\}$ denote real and imaginary parts, respectively, and $\underline{\varphi}$ is the complex-valued mode shape vector.

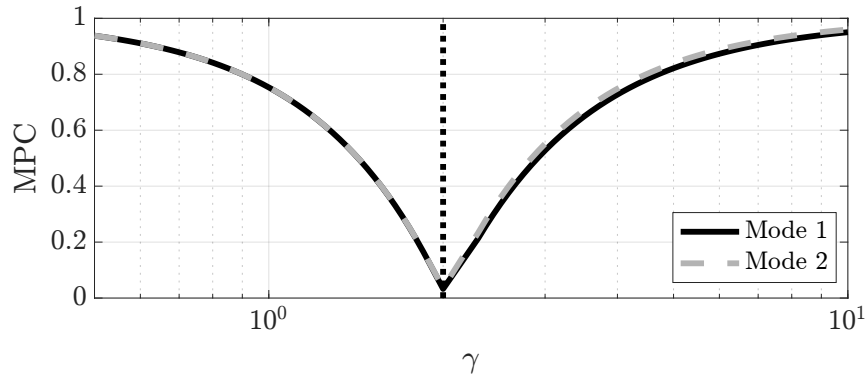

*Figure 2. Modal phase collinearity of vibration modes for varying $\gamma$ for the system with $\lambda = 0, \lambda_1 = 0.00125, \lambda_2 = 0.05$; the dashed line indicates the critical value $\gamma = 2$.*

Figure 2 shows the MPC results of the system at different $\gamma$ values; note that MPC $= 1$ corresponds to a real-valued mode shape vector. This indicates that the system vibration modes are most highly complex at $\gamma = 2$, and maintain a high degree of complexity around the neighborhood of $\gamma = 2$.

## 2.2. Impulsive responses, modal interactions and energy decay

Considering impulsive excitations applied to the system of Fig. 1 and zero initial conditions, its responses are expressed in state space as follows,

$$\begin{bmatrix} \underline{u}(t) \\ \underline{v}(t) \end{bmatrix} = [\Phi][\eta(t)] \quad (2)$$

where $[\Phi]$ is a 4x4 complex-valued mode shape matrix, $[\eta(t)]$ is the 4x1 complex-valued vector of modal responses, $\underline{u}(t)$ is the 2x1 displacement vector, and $\underline{v}(t)$ is the 2x1 velocity vector. Note that this representation expands the dynamic analysis to the four-dimensional state space. Then the elements of the modal vector are represented as follows,

$$\eta_k(t) = \alpha_k e^{-\sigma_k t + j\omega_k t}, k = 1, \dots, 4 \quad (3)$$

where $\alpha_k, \sigma_k, \omega_k$ are the k-th modal amplitude, dissipation rate, and frequency, respectively, and $j = \sqrt{-1}$ is the imaginary constant; for analytical expressions of these modal parameters depending on $\gamma$ being smaller or greater than 2, we refer to (Baldelomar Pinto et al., 2025) where a multiple-scales analysis for weak coupling and damping was performed. In synopsis, for $\gamma < 2$, there exist two distinct dissipation rates, $\sigma_{1,2} = \frac{\lambda_+}{2} \mp \frac{\lambda_-}{4}\sqrt{4-\gamma^2}$, with $\lambda_+ = \frac{\lambda_2 + \lambda_1}{2}, \lambda_- = \frac{\lambda_2 - \lambda_1}{2}$, and a single frequency $\omega_1 = \omega_2 = 1 + \frac{\beta}{2}$; whereas for $\gamma > 2$ there exists a single dissipation rate $\sigma_1 = \sigma_2 = \lambda_+$, and two distinct (but closely spaced) frequencies $\omega_{1,2} = \omega_0 \pm \omega_d$, where $\omega_d = \frac{\lambda_-}{4}\sqrt{\gamma^2 - 4}$ is the (slow) frequency governing the beat phenomena in the velocity envelopes. Moreover, the modal amplitudes $\alpha_k$ depend on the impulsive excitations that are applied at $t = 0+$ when the system is initially at rest.

Based on (3), we can construct expressions for the two oscillator velocities as follows,

$$v_i(t) = 2|\psi_{i1}|e^{-\sigma_1 t}\cos(\omega_1 t + \phi_{i1}) + 2|\psi_{i2}|e^{-\sigma_2 t}\cos(\omega_2 t + \phi_{i2}), \quad i = 1,2 \quad (4)$$

where $\psi_{im} = \Phi_{im}\alpha_m$, and $\phi_{im} = \arg\{\psi_{im}\}$. We note that up to this point we have not made any assumptions regarding damping non-proportionality in the system.

Expressions (4) represent the velocity responses of the two oscillators in terms of modal superpositions. Since we are interested in defining the velocity envelopes of these responses, we apply the Hilbert transform $\mathcal{H}[\bullet]$ in (4) to derive the following expressions,

$$\hat{v}_i(t) \equiv \mathcal{H}[v_i(t)] =$$
$$2|\psi_{i1}|e^{-\sigma_1 t}\sin(\omega_1 t + \phi_{i1}) + 2|\psi_{i2}|e^{-\sigma_2 t}\sin(\omega_2 t + \phi_{i2}), \quad i = 1,2 \tag{5}$$

which, in turn, we use to define analytic functions representing the velocity responses in the complex plane as, $z_i(t) \equiv v_i(t) + j\hat{v}_i(t)$:

$$z_i(t) = 2|\psi_{i1}|e^{-\sigma_1 t + j\omega_1 t + j\phi_{i1}} + 2|\psi_{i2}|e^{-\sigma_2 t + j\omega_2 t + j\phi_{i2}}, \quad i = 1,2 \tag{6}$$

Based on these complex velocity functions, we may define the envelope functions, $\langle v_i(t) \rangle$, of the corresponding velocities $v_i(t)$ by simply taking the moduli of $z_i(t)$:

$$\langle v_i(t) \rangle \equiv |z_i(t)| =$$
$$\sqrt{|\psi_{i1}|^2 e^{-2\sigma_1 t} + |\psi_{i2}|^2 e^{-2\sigma_2 t} + 2|\psi_{i1}||\psi_{i2}|e^{-(\sigma_1+\sigma_2)t}\cos((\omega_1-\omega_2)t + \phi_{i1} - \phi_{i2})}, i = 1,2 \tag{7}$$

where the $\langle \bullet \rangle$ operator indicates the envelope of the time series of the response.

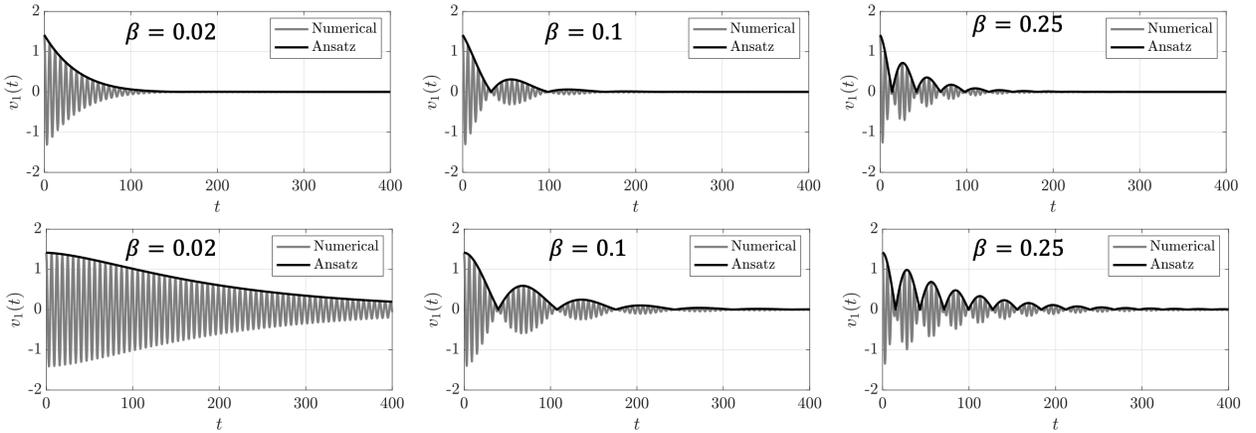

*Figure 3. Velocity response and corresponding velocity envelope of oscillator 1 subject to a single impulse of intensity $\sqrt{2}$ applied to it, for different coupling stiffnesses β: Top plots are for a classically (proportionally) damped system with $\lambda_1 = \lambda_2 = 0.05$, while bottom plots are for a non-classically (non-proportionally) damped system with $\lambda_1 = 0.00125, \lambda_2 = 0.05$.*

In Figure 3, we present the numerical verification of the envelope expression (7) for oscillator 1, when a single impulse of intensity $\sqrt{2}$ is applied to that oscillator. We note that the analytical envelopes (7) agree with the numerically computed envelopes for varying values of coupling stiffness $\beta$, irrespective of proportional or non-proportional viscous damping distribution.

It follows that (7) eliminates any ambiguity associated with approximating the envelope numerically (Mojahed et al., 2022a,b; Chang et al., 2025), and highlights the modal contributions to the envelopes; in turn, these expressions reveal the modal interactions that occur in the velocity responses due to the combined effects of closely spaced modes and damping non-proportionality in the neighborhood of $\gamma = 2$. In Fig. 3 the interactions between oscillators (in both rows of plots) and modes (only in the bottom row of plots) are shown as beat phenomena for relatively strong coupling – which correspond to $\gamma > 2$ (Baldelomar Pinto et al., 2025).

From the results of Fig. 3 we deduce that for relatively strong values of coupling $\beta$ (or, equivalently, for $\gamma > 2$) the velocity envelopes, as defined by fitting curves connecting the maxima of the velocity time series, cease from being smooth functions of time; this feature prevents the very definition of velocity envelope functions which, as reviewed below, are required for the computation of BW and EST of the two-DOF system of Fig. 1 close to the critical $\gamma = 2$ (Mojahed et al., 2022a,b; Chang et al., 2025). Motivated by these results, we replace these velocity envelopes with the *total energy decay* in the two-DOF system, which is a smooth, well-defined and monotonically decaying function (since the system is passive, non-self-excited, and only dissipates energy in time). Moreover, this total energy measure not only accounts for the possible occurrence of modal interactions due to non-classical damping distribution, but also for different types of impulse excitation, i.e., excitation of only oscillator 1 or 2, or of both oscillators simultaneously. With this in consideration, we analytically approximate the total energy decay in the system using the velocity expressions (4), and then obtain the following total energy decay

$$E(t) = 2(|\psi_{11}|^2 + |\psi_{21}|^2)e^{-2\sigma_1 t} + 2(|\psi_{12}|^2 + |\psi_{22}|^2)e^{-2\sigma_2 t} +$$
$$4e^{-(\sigma_1+\sigma_2)t} \left[\, |\psi_{11}||\psi_{12}| \cos\big((\omega_2 - \omega_1)t + \phi_{11} - \phi_{12}\big) + \right.$$
$$\left. |\psi_{21}||\psi_{22}| \cos\big((\omega_2 - \omega_1)t + \phi_{21} - \phi_{22}\big) \,\right] \quad (8)$$

In Figure 4, we verify the validity of the total energy decay expression (8) for the non-classically damped system with $\gamma$ less or greater than the critical value 2. The two distinct regimes of the decaying dynamics following the initial application of the impulse are clearly deduced in these results. From here on, we will distinguish between two main types of impulsive excitation: Case 1 refers to a single impulsive excitation, $\sqrt{2}\delta(t)$, applied to the lightly damped oscillator 1 with the system being initially at rest, whereas Case 2 corresponds to the same impulsive excitation, $\sqrt{2}\delta(t)$, applied only to the heavily damped oscillator 2.

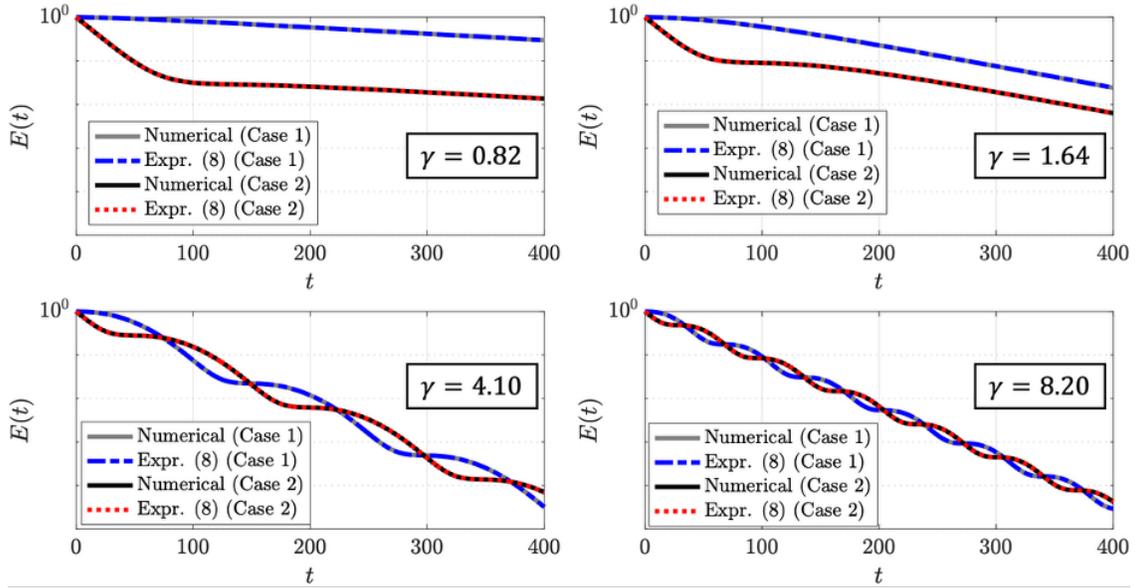

*Figure 4. Validation of the energy expression (8) for the two-DOF system with non-classical damping and closely spaced modes, for varying γ and two types of impulsive excitation.*

From the results of Fig. 4, the analytical expression (4) is validated with direct numerical simulations. Focusing on the individual plots of Fig. 4, for $\gamma < 2$ the existence of two distinct dissipation rates (each dominating early and later time decay) and a single oscillation frequency, is apparent; the total energy decay is a smooth monotonically decaying function of time with no oscillations. The dynamics change drastically for $\gamma > 2$, where a single dissipation rate and two distinct (but closely spaced) frequencies exist; then the total energy decay is oscillatory due to energy exchanges between oscillators (and interacting modes). But even in this latter case, the total energy decay is a smoothly decaying function which, as shown below, can be used as basis for investigating the overall dissipative capacity (in terms of BW) and dissipation rate (in terms of EST) of the two-DOF system of Fig. 1. We conclude by mentioning that for a proportionally damped system, e.g., setting $\lambda = 0, \lambda_1 = \lambda_2$ in Fig. 1, the total energy is exponentially decaying, being dependent on the energy decay of each individual mode (which cannot interact and exchange energy with any other mode).

## 3. Effective oscillator, bandwidth (BW) and energy dissipation rate (EST)

### 3.1. Definitions and results

Recent work (Mojahed *et al.*, 2022a) generalized the classical concept of BW defined for a single-DOF LTI resonator (Ewins, 2009), to a more general class of single-DOF nonlinear or time-variant resonators. To achieve this, previous signal processing results (Gabor, 1949; Starosielec and Hägele, 2014) were extended to dynamics by defining the BW of a general single-DOF resonator as,

$$\Delta \omega = 2\sqrt{\frac{\int_0^\infty (\omega - \omega_c)^2 |V(\omega)|^4 d\omega}{\int_0^\infty |V(\omega)|^4 d\omega}} \qquad (9)$$

where, $V(\omega) = \mathcal{F}[v(t)]$ is the Fourier transform of the velocity envelope, $v(t)$, of the impulsive response of the resonator, and $\omega_c$ is the so-called central frequency where the maximum of $|V(\omega)|$ occurs, with $\omega_c = 0$ for a monotonically decreasing envelope $v(t)$. Note that $|V(\omega)|^2$ is proportional to the power dissipated by the resonator, so (9) provides a measure of energy dispersion in the frequency domain, or a measure of dissipative capacity, that is, the frequency range where the energy of the resonator is most effectively dissipated. For a single-DOF LTI resonator expression (3) recovers the classical half-power BW.

An analogous expression can be derived for a general single-DOF resonator, to characterize the locality of the square of the envelope of the impulsive velocity response in the time domain, which provides the energy storage time (EST), or a measure of the rate of energy dissipation of the resonator in the time domain:

$$\Delta t = \sqrt{2} \sqrt{\frac{\int_0^\infty v(t)^4 t^2 dt}{\int_0^\infty v(t)^4 dt}} \qquad (10)$$

Again, for a single-DOF LTI resonator, expression (10) recovers the classical time constant, $\Delta t = 1/\Delta \omega$, so that the fundamental TBL $\Delta t \Delta \omega = 1$ holds for this system.

Based on these definitions, we represent the two-DOF system of Fig. 1 as an *effective single-DOF oscillator* whose energy decay is equal to the total energy decay (8) of the two-DOF system. Clearly, this effective oscillator is not time-invariant (it is a fictitious oscillator after all that accounts for the combined energy dissipation by the two vibration modes of the system), yet we can define its *effective velocity envelope* as,

$$\langle v_{eff}(t) \rangle \equiv \left(\sqrt{2/M}\right)\sqrt{\langle E(t) \rangle} \qquad (11)$$

where $M$ is the total mass of the two-DOF system. The rationale is that the decay envelope of the total kinetic energy of the effective oscillator is identical to the decay envelope of its total energy, which, in turn, is set equal to (8). But then, using the definition of effective velocity envelope (11), we can proceed to compute the BW and EST of the effective oscillator through relations (9) and (10), respectively, and study the dissipative capacity and dissipation rate of the two-DOF system when viewed as a whole, i.e., in terms of its overall capacity to dissipate and store energy in the frequency and time domains. We note that this methodology can be extended to a more general class of classically or non-classically damped, multi-DOF oscillators, that include nonlinear and time-varying effects. Also, through (11) modal interactions are fully considered, in the sense that their effects regarding the overall energy dissipation capacity of the system are accounted for. Moreover, it is clear that since the effective oscillator is not time-invariant, its time bandwidth product (TBP) will not satisfy the fundamental TBL.

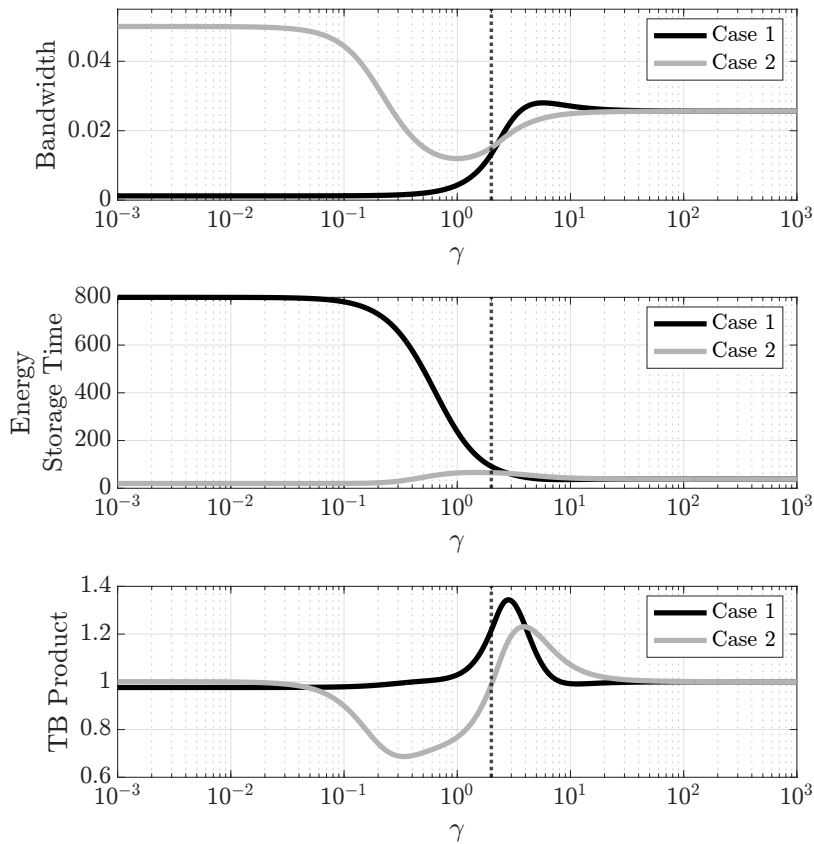

*Figure 5. BW, EST and TBP for varying γ of the effective oscillator of the two-DOF non-classically damped system subject to impulsive excitation of oscillator 1 – Case 1, or oscillator 2 – Case 2 (non-dimensional frequency and time units are used due to normalized parameters).*

Exceeding the TBL from above or below for the effective oscillator corresponding to the two-DOF system of Fig. 1 has interesting implications regarding its overall dissipative capacity and energy storage time. Indeed, if it would hold that $\Delta t\, \Delta \omega > 1$, we could achieve larger-bandwidth resonance – equivalently, larger dissipative capacity with the two-DOF system compared to a single-DOF LTI resonator, while maintaining the same EST (time constant); or, conversely, smaller EST compared to a single-DOF LTI resonator with the same BW while maintaining the same dissipative capacity. The reverse effects would occur for resonators with $\Delta t\, \Delta\omega < 1$. Such comparisons are performed in a later section of this work.

In Fig. 5 we present the numerically computed BW, EST and TBP of the effective oscillator corresponding to the two-DOF non-classically damped system when a single impulsive excitation of intensity $\sqrt{2}$ is applied either to the lightly damped oscillator 1 (Case 1) or the heavily damped oscillator 2 (Case 2). In these simulations we systematically vary the coupling stiffness values, $\beta$, while maintaining constant damping values, $\lambda = 0, \lambda_1 = 0.00125$ and $\lambda_2 = 0.05$, ensuring non-classical damping distribution and effectively changing $\gamma$. At each $\gamma$ value, we compute the eigenvalues of the system, gather the modal parameters to determine the total energy envelope (8), and based on this, define the effective oscillator velocity envelope through (11). Following this, we numerically compute the BW and EST and the TBP according to expressions (9) and (10).

The results of Fig. 5 reveal the sensitivity of the overall dissipative capacity and rate of energy dissipation of the two-DOF non-classically damped system at the point of application of the impulsive excitation. Indeed, when the impulse is applied to the right, more heavily damped oscillator (see Fig. 1 – Case 2), as $\gamma$ increases from low to high values, the overall BW of the system ranges from $\Delta\omega \to 0.05 \equiv \lambda_2$ for $\gamma \ll 1$ to $\Delta\omega \to 0.0256 \equiv \frac{\lambda_1+\lambda_2}{2}$ for $\gamma \gg 1$; a reverse variation occurs when the impulse is applied to the more lightly damped oscillator (Case 1), where the BW ranges from $\Delta\omega \to 0.00125 \equiv \lambda_1$ for $\gamma \ll 1$ to $\Delta\omega \to 0.0256 \equiv \frac{\lambda_1+\lambda_2}{2}$ for $\gamma \gg 1$. This makes physical sense since for $\gamma \ll 1$ the coupling between the two oscillators vanishes and the two-DOF system of Fig. 1 degenerates to two uncoupled single-DOF oscillators, each with its own damping coefficient $\lambda_1$ or $\lambda_2$; this also explains why the TBP reaches the fundamental TBL $\Delta t \Delta\omega = 1$ for $\gamma \ll 1$. Similarly, for $\gamma \gg 1$ the coupling between oscillators becomes nearly rigid, the two-DOF degenerates to a single-DOF oscillator albeit with damping coefficient equal to $\frac{\lambda_1+\lambda_2}{2}$, and again exhibits the fundamental TBL $\Delta t \Delta\omega = 1$.

Similar, but inverse patterns occur for the EST results. That is, for low $\gamma$ the EST of the effective oscillator tends to the inverse of the damping coefficient associated with the directly excited oscillator, $\Delta t \to 1/\lambda_i$, whereas, for high $\gamma$ the EST tends to the inverse of the average of the two damping coefficients, $\Delta t \to \frac{2}{\lambda_1+\lambda_2}$.

It follows from these results that *the combined effects of closely spaced modes and non-classical damping distribution are realized mainly in the neighborhood of $\gamma = 2$, where the TBP $\Delta t \Delta \omega$ is exceeded from above for Case 1 and from both above and below for Case 2*. Furthermore, the results highlight the fact that the only dynamical regime where the system of Fig. 1 behaves as a "genuine" two-dimensional dynamical system is in the neighborhood of $\gamma = 2$, with everywhere else the system exhibits predominantly a set of two nearly uncoupled single-DOF oscillators (for $\gamma \ll 1$) or a single-DOF oscillator (for $\gamma \gg 1$). In addition, we note that in the neighborhood of $\gamma = 2$, the two vibration modes have the highest degree of complexity (see Fig. 2), and this is precisely the dynamical regime where the TBP deviates the most from unity (the fundamental TBL $\Delta t \Delta \omega = 1$). This indicates that mode complexity due to the combined effects of mode proximity and damping non-proportionality affects greatly the overall dissipative capacity of the two-DOF system and its rate of energy dissipation (or energy storage capacity) in time, resulting in a significant deviation of the dynamics from single-DOF behavior.

### 3.2. Transition towards classical damping distribution

To elucidate the transition to a proportional damping distribution in the two-DOF system and the resulting effects on the BW $\Delta \omega$, EST $\Delta t$ and TBP $\Delta t \Delta \omega$, we performed a numerical study by fixing the damping coefficient $\lambda_2 = 0.05$ of oscillator 2, incrementally varying the damping of oscillator 1, $\lambda_1$, in the range $[0.00125 - 0.05]$, and continuously varying the coupling stiffness $\beta$ (thus, in effect varying $\gamma$). Clearly, in the upper limit of the range of variation of $\lambda_1$ the system becomes proportionally damped, so by this exercise we study the effect of non-proportional (non-classical) viscous damping distribution on the dissipative properties of the system, and most importantly on the TBP.

Next, we compute the effective velocity envelope $\langle v_{eff}(t) \rangle$ through (11) for different values of $\gamma$, and based on that, the corresponding BW, EST and TBP. The results, shown in Figure 6, indicate that by increasing $\lambda_1$ to its upper bound the BW eventually converges to the value of both damping coefficients, $\lambda_1 = \lambda_2 = 0.05$, whereas the EST converges to the inverse of the damping coefficients, $1/\lambda_1 = 1/\lambda_2 = 20$. It follows that in the limit of proportional damping distribution the TBP of the two-DOF system convergences to unity, which is the fundamental TBL for single-DOF LTI resonators. We note that this limiting behavior is realized regardless of initial conditions, which is expected since both oscillators possess identical damping coefficients. Hence, the two modes are uncoupled and the two-DOF system dissipates energy similar to a SDOF LTI resonator.

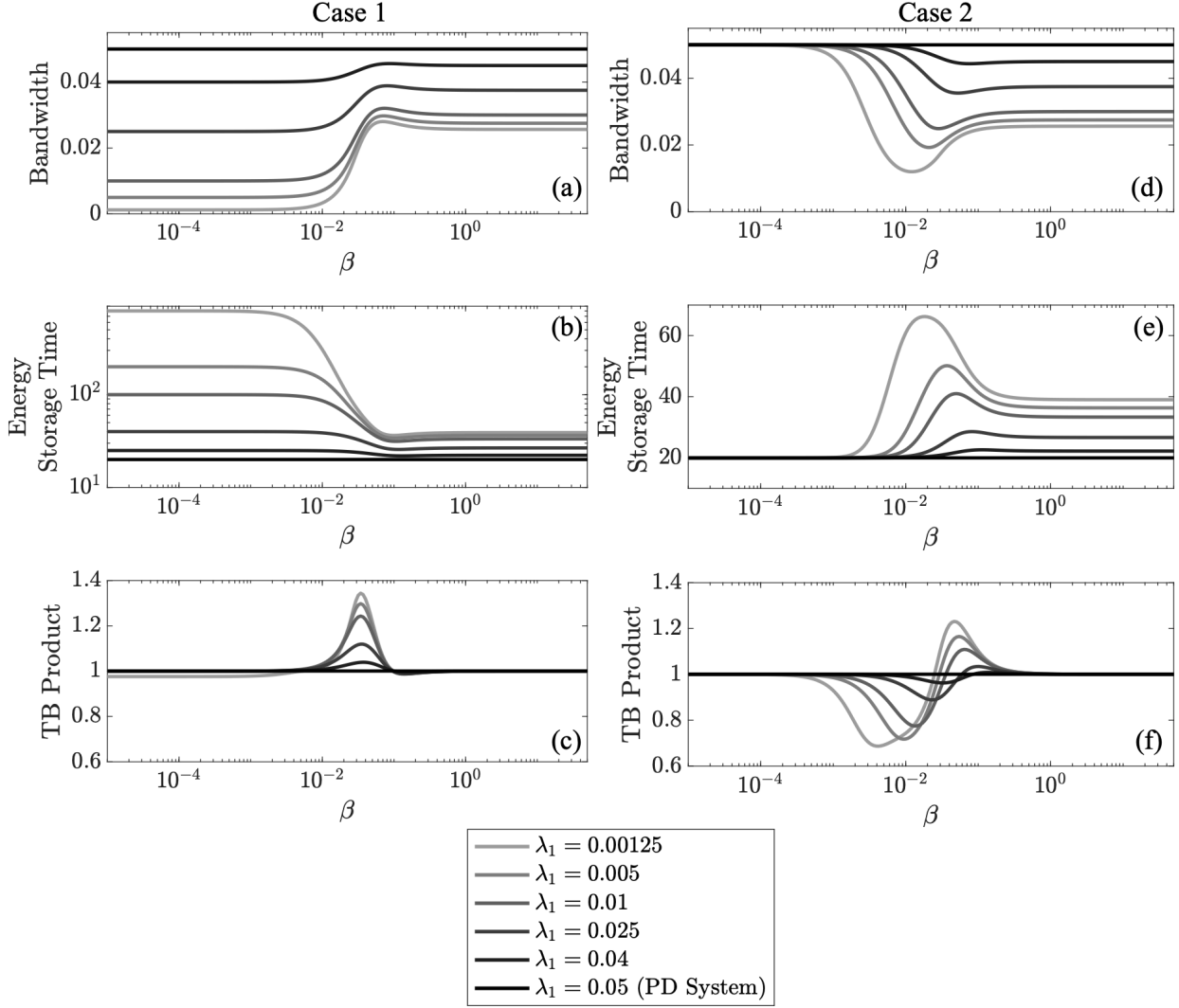

Figure 6. Transition from non-proportional to proportional damping distribution in the two-DOF system of Fig. 1: BW, EST and TBP of the effective oscillator with fixed $\lambda_2 = 0.05$, and varying coupling $\beta$ and damping coefficient $\lambda_1$, for Case 1 (a-c) and Case 2 (d-f) of impulsive excitation.

### 3.3. Extreme break of the fundamental TBL in the two-DOF system

To gain physical insight into the transient responses of the two-DOF system at values of $\gamma$ corresponding to extreme break of the fundamental TBL, $\Delta t \Delta \omega = 1$, we numerically investigate the transient dynamics of the system energy response of the 2DOF system subject to an impulsive excitation applied to either oscillator 1 (Case 1) or oscillator 2 (Case 2) for two specific values of $\gamma$. Specifically, referring to the plots of Fig. 5, we select the value $\gamma = 0.35$ corresponding to the *minimum of the TBP, $\Delta t \Delta \omega = 0.687 < 1$ for Case 2, and the value $\gamma = 2.70$ corresponding to the maximum of the TBP, $\Delta t \Delta \omega = 1.341 > 1$ for Case 1.*

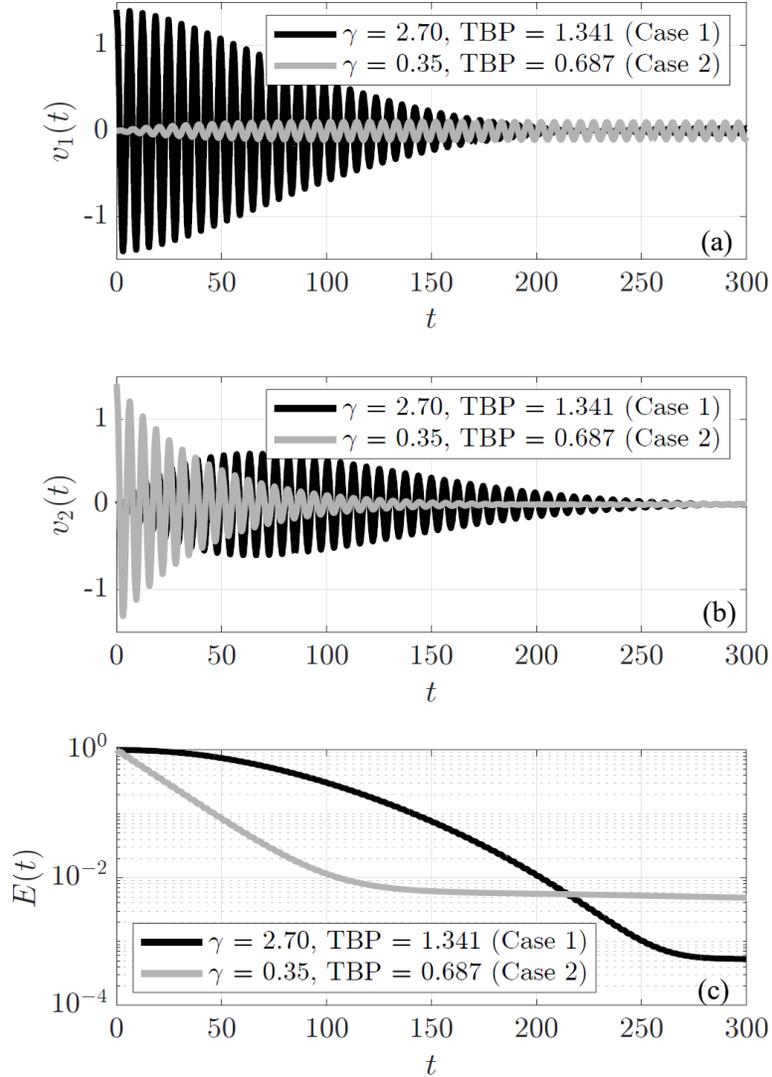

*Figure 7. Impulsive responses of the non-classically damped two-DOF for $\gamma = 0.35$ corresponding to minimum TBP for Case 2 (gray line) and $\gamma = 2.70$ corresponding to maximum TBP for Case 1 (black line): (a) Velocity of oscillator 1, (b) velocity of oscillator 2, and (c) total energy decay in the two-DOF system.*

From the results we deduce that the two-DOF system with $\Delta t \Delta \omega < 1$ dissipates energy much more efficiently than the system with $\Delta t \Delta \omega > 1$ – compare the response of the impulsively excited oscillator 1 for Case 1 (Fig. 7a) to the response of the impulsively excited oscillator 2 for Case 2 (Fig. 7b). Moreover, this is more clearly discerned when comparing the total energy decay plots in Fig. 7c, when the two-DOF system with TBP $\Delta t \Delta \omega = 0.687$ (Case 2) dissipates energy more efficiently compared to the system with TBP $\Delta t \Delta \omega = 1.341$ (Case 1). Conversely, the two-DOF system with $\Delta t \Delta \omega > 1$ stores energy for longer duration compared to the system with $\Delta t \Delta \omega < 1$. These results indicate that by using the energy dissipation measures developed in this work one can predictvely design a general multi-DOF system (even with nonlinear and/or time-varying terms) to dissipate energy more (or less) efficiently and/or at a faster (or slower) rate, compared to a linear SDOF LTI resonator with the same mass and identical BW or EST. This comparison is further studied in the next section.

### 3.4. Comparison with a single-DOF LTI resonator

The previous results concerning the break (from above or below) of the fundamental TBL for the two-DOF non-classically damped system dictate a comparative study of its dissipative capacity and rate of energy dissipation with a single-DOF LTI resonator which always satisfies the fundamental TBL, $\Delta t \Delta \omega = 1$, regardless of mass, stiffness, and damping coefficients. Accordingly, as a first step we examine the total energy decay in the two-DOF system and compare it to the corresponding energy decay in an underdamped single-DOF resonator (Fig. 8) with normalized parameters $m = \beta = 1$ and normalized damping coefficient $\lambda$ being equal to the corresponding bandwidth $\Delta \omega$ of the two-DOF system. For this comparison we select the same $\gamma$ values as in the previous section, which correspond to the maximum (for Case 1) and minimum (for Case 2) TBP of Fig. 5.

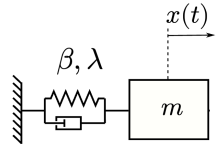

*Figure 8. Single-DOF LTI oscillator used in comparative study.*

In Fig. 9 we compare the total energy decays in the two-DOF system and the same-BW single-DOF LTI resonator. For $\gamma = 2.7$ with $\Delta t \Delta \omega = 1.341 > 1$ (Case 1, Fig. 9a) the total energy of the two-DOF system decreases more slowly ($\Delta t = 65.4$) compared to the single-DOF LTI oscillator with the same bandwidth ($\Delta t = 48.78$). For this value of $\gamma$ there exist intense beat phenomena, yielding energy exchanges between the two oscillators and modes, which delay the overall energy dissipation. For $\gamma = 0.35$ with $\Delta t \Delta \omega = 0.687 < 1$ (Case 2, Fig. 9b) we observe that the energy of the two-DOF system is dissipated more rapidly ($\Delta t = 33.68$) than in the single-DOF LTI system with the same bandwidth ($\Delta t = 49.02$). This highlights how the existence of two distinct dissipation rates in the two-DOF system contributes to faster energy dissipation than the single-DOF LTI system with a single dissipation rate. From these trends, we can determine that the TBP $\Delta t \Delta \omega$ serves as a useful measure for designing systems or devices with prescribed dissipative capacities and rates of energy dissipation, e.g., with capacity to dissipate energy faster (for $\Delta t \Delta \omega < 1$) or store energy for longer times (for $\Delta t \Delta \omega > 1$) compared to the predictions of the fundamental TBL which is only valid for single-DOF LTI resonators.

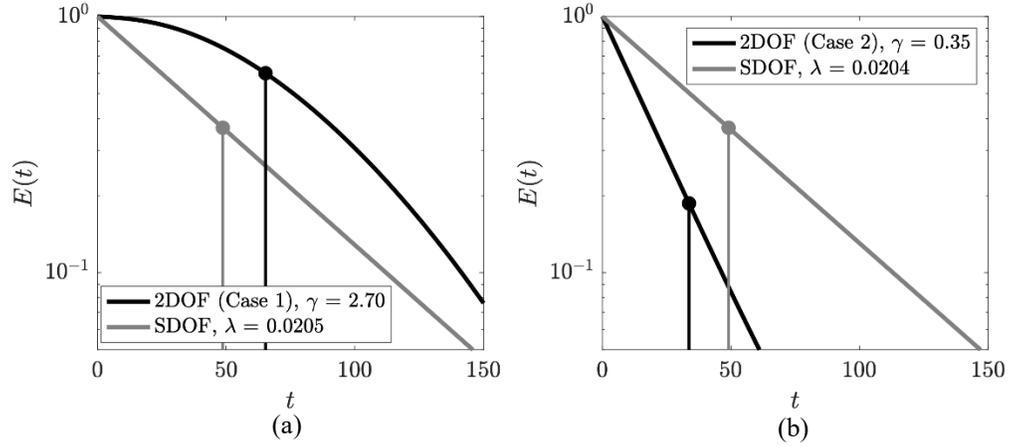

Figure 9. Comparison of total energy decays of the two-DOF system (black line) and a single-DOF oscillator of the same BW $\Delta\omega$ (gray line): (a) Case 1 for $\gamma = 2.70$ and TBP $\Delta t \Delta \omega > 1$, and (b) Case 2 for $\gamma = 0.35$ with TBP $\Delta t \Delta \omega < 1$ (vertical lines indicate the corresponding EST $\Delta t$).

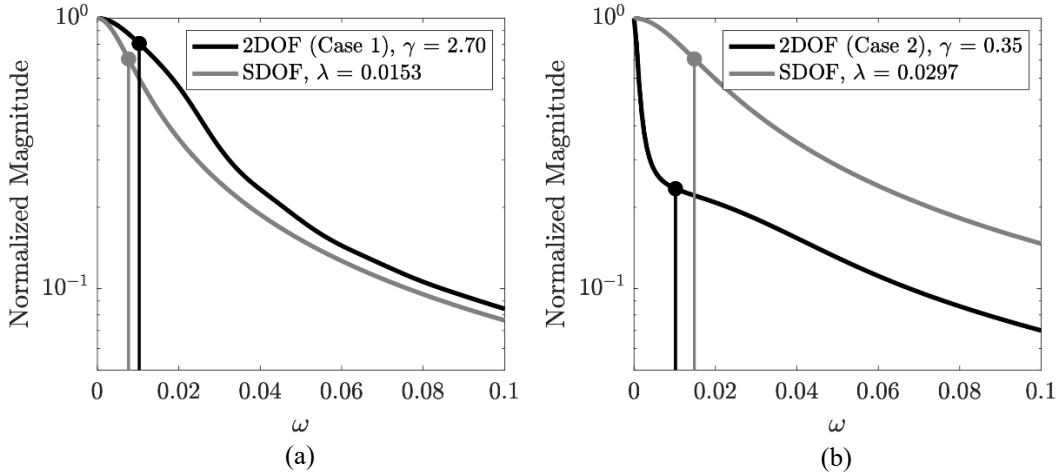

Figure 10. Comparison of the moduli of the FFTs of the velocity envelopes of the two-DOF system (black line) and a single-DOF oscillator of the same EST $\Delta t$ (gray line): (a) Case 1 for $\gamma = 2.70$ and TBP $\Delta t \Delta \omega > 1$, and (b) Case 2 for $\gamma = 0.35$ with TBP $\Delta t \Delta \omega < 1$ (vertical lines indicate the corresponding half-BW $\Delta\omega/2$).

Similarly, we compare the BW of the two-DOF system with that of a single-DOF LTI resonator with the same EST $\Delta t$ (with parameters $m = \beta = 1, \lambda = 1/\Delta t$). To this end, in Fig. 10 we compare the moduli of the FFT of the velocity envelopes of the two systems and indicate the corresponding BW. For $\gamma = 2.7$ with $\Delta t \Delta \omega = 1.341 > 1$ (Case 1, Fig. 10a) we observe that the frequency spectrum for the two-DOF system has a wider spread about $\omega = 0$ ($\Delta \omega = 0.0206$), i.e., is more broadband, compared to the single-DOF LTI resonator with the same time constant ($\Delta \omega = 0.0153$). On the contrary, for $\gamma = 0.35$ with $\Delta t \Delta \omega = 0.687 < 1$ (Case 2, Fig. 10b), the single-DOF LTI resonator with identical time constant has a wider spread in its frequency spectrum ($\Delta \omega = 0.0297$) compared to the more narrowband two-DOF system ($\Delta \omega = 0.02$). It follows that for $\Delta t \Delta \omega > 1$ ($< 1$) the two-DOF is more broadband (narrowband) than the single-DOF LTI resonator with the same dissipation rate or EST. These results pave the way for designing, e.g., broadband systems or devices that, while maintaining a set energy dissipation rate, are more broadband or narrowband compared to the prediction of the fundamental TBL.

## 4. Experimental Validation

Experimental tests were performed to validate the theoretical predictions. Fig. 11 shows the experimental setup which is the same used in Baldelomar Pinto et al. (2025) to study the distinct dynamical regimes in a realization of the two-DOF system of Fig. 1. In that work, modal analysis was applied to construct a reduced-order model of the experimental fixture, which is similar to Fig. 1, corresponding to different realizations of the coupling to damping non-proportionality parameter $\gamma$.

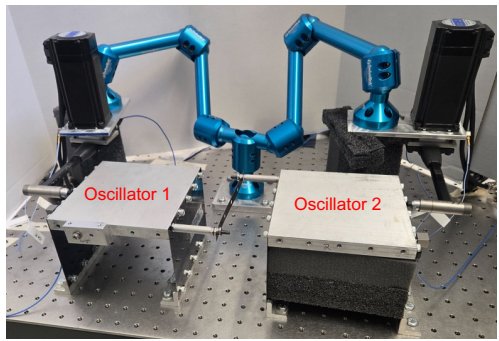

*Figure 11. The experimental system consisting of a pair of near-identical harmonic oscillators (but for the damping) coupled with a thin steel flexure and excited by a pair of automatic modal hammers (Baldelomar Pinto et al., 2025).*

Using experimental modal analysis, a reduced order model (ROM) for the experimental setup was obtained having the configuration of Fig. 1, and the following system parameters (Baldelomar Pinto et al., 2025). Mass of oscillator 1 was $m_1 = 1.073\ Kg$ and of oscillator 2, $m_2 = 1.094\ Kg$. Stiffness of oscillator 1 was $k_1 = 30019$ N/m and grounding damping coefficient $c_1 = 0.614$ Ns/m, while $k_2 = 30606$ N/m and $c_2 = 16.51$ Ns/m, respectively, for oscillator 2 (yielding non-dimensional stiffness and damping parameters for oscillator 1, $\beta_1 = 1$ and $\lambda_1 = 0.003$, and $\beta_2 = 1.02$, $\lambda_2 = 0.092$, respectively, for oscillator 2). Lastly, two values of coupling stiffness were used, namely, $K = 1181.4$ N/m and $K = 2503$ N/m, respectively, yielding non-dimensional coupling stiffnesses $\beta = 0.039$ and $\beta = 0.083$. Based on these parameters two different values of the $\gamma$ parameter are experimentally realized, namely, $\gamma = 1.78$ and $\gamma = 9.27$.

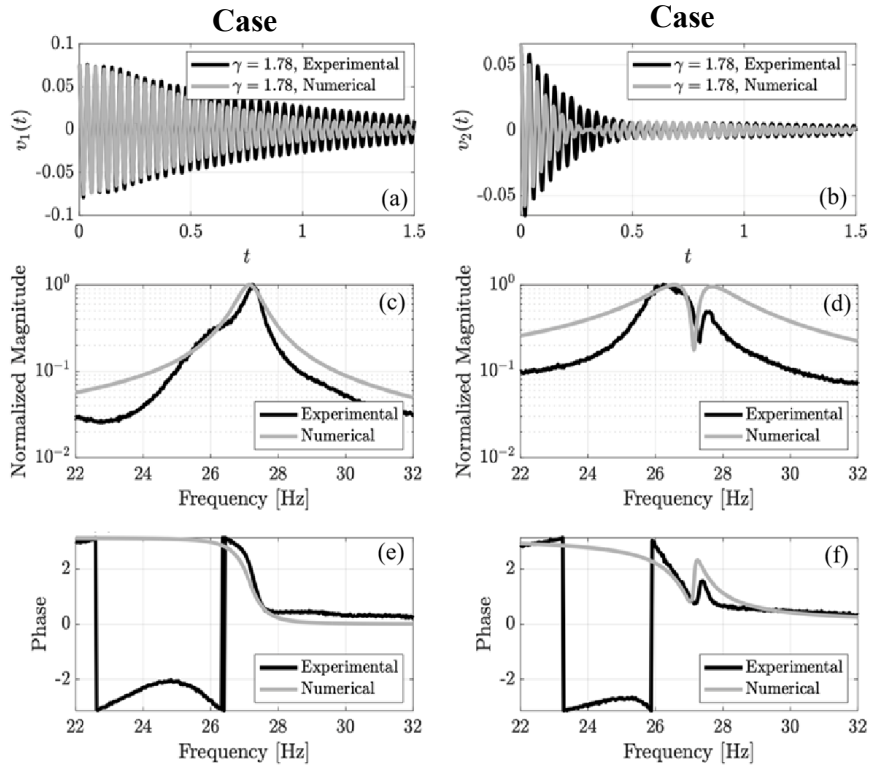

*Figure 12. Experimental (black line) and ROM (gray line) impulsive responses from the fixture with γ = 1.78: Transient velocities for (a) Case 1, and (b) Case 2; modulus of receptance for (c) Case 1 and (d) Case 2; and phase of receptance for (e) Case 1 and (f) Case 2.*

In Figs. 12 and 13 we depict comparisons between the experimental impulsive transient and numerical simulations of the corresponding ROMs for the fixture with weaker and stronger coupling corresponding to $\gamma = 1.78$ and $\gamma = 9.27$, respectively. Based on the observed agreement between experimental and numerical results, the validity of the experimentally identified ROM is established. For these simulations, the initial displacements of the ROM were set equal to zero, and the initial velocities (equivalently, the initially applied impulses) were set equal to the maximum velocity amplitudes measured at the corresponding experiments; this ensured proper comparison between the ROM simulations and the experimental measurements. The differences between the experiment and ROM results can be attributed to unmodeled or uncertainty effects, e.g., errors in the identified modal parameters and especially the damping estimation, and modeling uncertainty in the damping of the coupling element of the fixture which in the ROM is taken to be equal to zero.

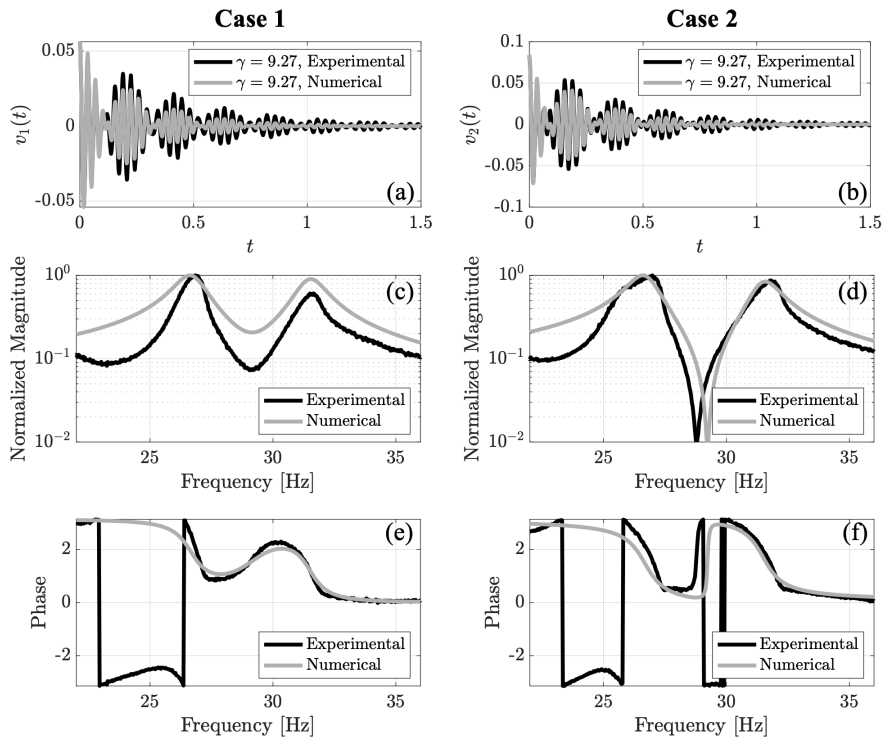

*Figure 13. Experimental (black line) and ROM (gray line) impulsive responses from the fixture with $\gamma = 9.27$: Transient velocities for (a) Case 1, and (b) Case 2; modulus of receptance for (c) Case 1 and (d) Case 2; and phase of receptance for (e) Case 1 and (f) Case 2.*

Considering the velocity measurements of Figs. 12 and 13, we follow two approaches for extracting the BW and EST from each impulsive response. The first approach ("ROM approach") involves solving the eigenvalue problem associated with each ROM to extract modal parameters from each case and compute the total energy decay based on expression (8). Then, the corresponding effective velocity envelope is computed through (11) and used to determine the BW and EST with relations (9) and (10). The second approach ("the experimental method") involves performing numerical enveloping of the velocity time series of each oscillator by fitting the local maxima of its velocity time series with Akima splines. The velocity envelopes obtained from this process are then used to construct the experimental total energy decay of the system as, $E(t) = \frac{1}{2}m_1\langle v_1(t)\rangle^2 + \frac{1}{2}m_2\langle v_2(t)\rangle^2$, from which an effective velocity envelope is extracted and used to estimate BW and EST.

Table 1 presents the results for the BW, EST and TBP of the systems with $\gamma = 1.78$ (Fig. 12) and $\gamma = 9.27$ (Fig. 13). We note that the estimates for the BW and EST derived by the two approaches differ, which can be attributed to the fact that the system parameters of the ROM are obtained from experimental modal analysis, and, hence, are expected to contain errors in their estimation. However, it is noteworthy that resulting TBP values are nearly identical for the same $\gamma$ value and subject to the same initial conditions. It follows that the experimental estimation of the TBP is rather robust, and, hence, provides a reliable indicator for characterizing the overall dissipative and energy storage capacity of the system of coupled oscillators.

Table 1. BW, EST and TBP for the impulsive tests of Figs. 12 and 13.

| Case | $\gamma$ | ROM Approach | | | Experimental Approach | | |
|---|---|---|---|---|---|---|---|
| | | $\Delta\omega$ [rad/s] | $\Delta t$ [s] | TBP [-] | $\Delta\omega$ [rad/s] | $\Delta t$ [s] | TBP [-] |
| 1 | 1.78 | 3.4134 | 0.339 | 1.157 | 2.1884 | 0.491 | 1.075 |
| 1 | 9.27 | 8.3568 | 0.119 | 0.993 | 5.0722 | 0.196 | 0.994 |
| 2 | 1.78 | 4.4811 | 0.205 | 0.918 | 3.7982 | 0.216 | 0.819 |
| 2 | 9.27 | 7.5536 | 0.144 | 1.089 | 5.4010 | 0.186 | 1.005 |

Lastly, in Fig.14 we depict the (normalized with respect to initial energy) total energy decay of the ROM using (8) and of the experimental system fitting the measured velocity envelopes and using the relation $E(t) = \frac{1}{2}m_1\langle v_1(t)\rangle^2 + \frac{1}{2}m_2\langle v_2(t)\rangle^2$. From the results, we note good qualitative and quantitative agreement between the ROM and experimental data (compare Figs 14a and c, and Fig. 14b and d). Moreover, consistent with the theoretical predictions of the previous sections, faster energy dissipation occurs for the cases with $\Delta\omega\,\Delta t < 1$, while longer energy storage is observed for the cases with $\Delta\omega \Delta t > 1$. Indeed, the trends of total energy decay match those for the idealized two-DOF system of Fig. 1 for the same $\gamma$ values – e.g., compare the results of Figs. 14b,d with the total energy decay plots of Fig. 7c.

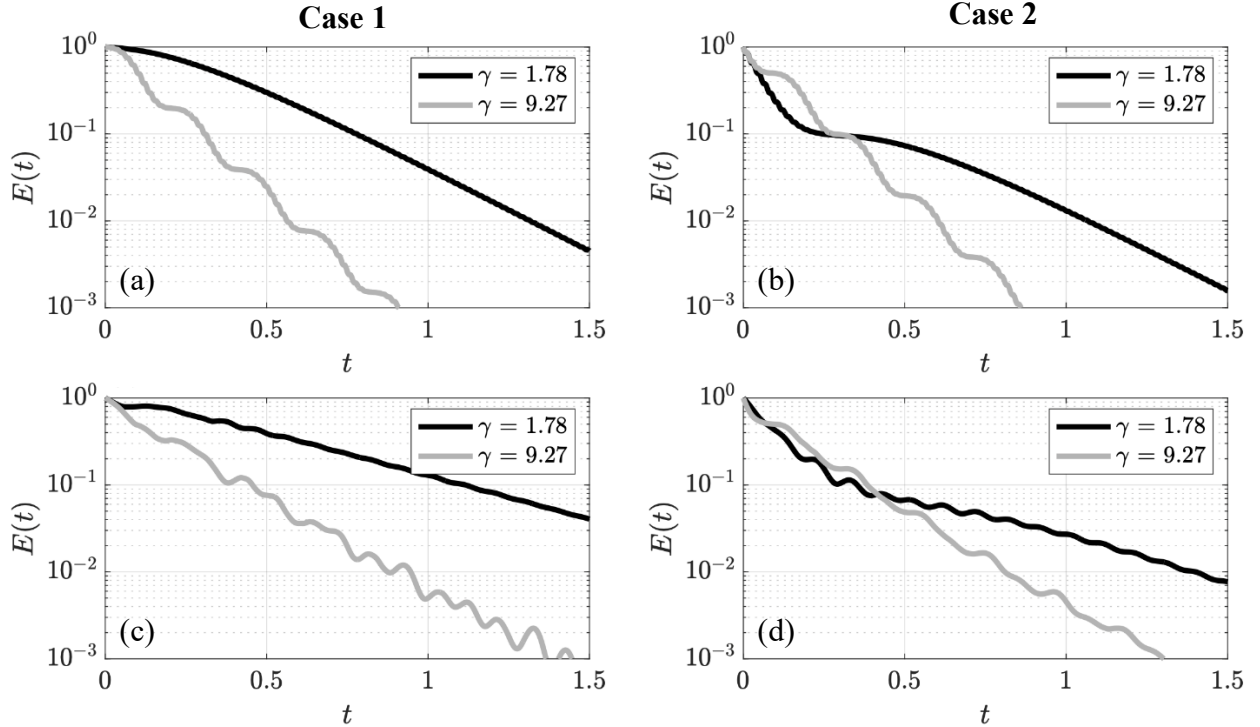

Figure 14. Normalized total energy decay from ROM (a,b) and experiment (c,d) for $\gamma = 1.78$ (black line) and $\gamma = 9.27$ (gray line); plots on the left correspond to Case 1 of impulse excitation, and on the right to Case 2.

## 5. Concluding remarks

We proposed scalar measures to quantify the dissipative capacity and rate of energy dissipation in an impulsively forced two-DOF non-classically damped system with closely spaced modes. The main challenge in the dynamics of this system is the occurence of modal interactions, that is, the coupling of the two vibration modes caused by the combined effects of non-proportional damping distribution and mode proximity. Apart from the fact that these modal interactions introduce two distinct regimes in the impulsive dynamics, i.e., a regime where a single (slow) energy modal exchange is realized, and another where recurrent modal energy exchanges take place (Baldelomar Pinto et al., 2025), modal interactions give rise to non-smooth velocity envelopes, and, hence, render inapplicable methodologies developed previously for single-DOF oscillators (Mojahed et al., 2022a) or multi-DOF proportionally damped systems (Chang et al., 2025). To address this issue we resort to defining an "effective velocity envelope" based on the total energy decay in the system, which being a smooth, monotonically decaying function of time can be used to compute effective bandwidth (BW) $\Delta\omega$, energy storage time $\Delta t$ (EST) and time-bandwidth product (TBP) $\Delta\omega\Delta t$.

Treating the two-DOF system as an effective oscillator, we studied its BW and EST for two cases of impulsive excitation and varying coupling stiffnesses. We showed that in the dynamical regime where modal interactions occur due to the combined effects of non-classical damping distribution and mode proximity, the impulsive response of the system exceeds the fundamental time-bandwidth limit (TBL) (being equal to unity), either from above or from below. This was to be expected since the TBL applies only to single-DOF LTI, lightly damped oscillators. Away from the regime of modal interactions, however, the TBP of the system approaches the TBL, and, for fixed damping, the two-DOF system behaves approximately as a set of two uncoupled single-DOF oscillators (for weak coupling), or as a single-DOF oscillator (for strong coupling).

In addition, the measures derived herein were used to assess the dissipation capacity of the two-DOF system. It was shown that two-DOF systems with $\Delta\omega\Delta t > 1$, i.e., exceeding the TBL from above, can store vibration energy for longer times compared to a single-DOF LTI resonator with the same BW; or conversely, they are more broadband than a single-DOF LTI resonator with the same energy dissipation rate. On the contrary, systems with $\Delta\omega\Delta t < 1$ dissipate energy much more efficiently and rapidly compared to a single-DOF LTI oscillator with the same BW; or, they are more narrowband while maintaining the same rate of dissipation with a single-DOF LTI resonator.

The theoretical predictions were verified experimentally and indicated that the TBP is a dissipative measure of practical significance. Moreover, we believe the results reported in this work shed insight into the mechanisms of energy dissipation in non-classically damped systems with modal interactions. Given that the methodology outlined in this work applies for general classes of linear or nonlinear, single- or multi-DOF systems, and time varying or time invariant parameters, the quantitative scalar measures developed here can be used to assess and optimize the overall inherent dissipative capacity of a dynamical system for changes of its parameters. Clearly, though, in the presence of nonlinear or time-varying effects, apart from the dependence of the dissipative measures on the location or locations of the applied impulsive excitations, these will depend also on energy. This calls for a new framework for quantifying and optimizing the dissipation properties of such general classes of dynamical systems.